\documentclass{article}

\usepackage{amsmath, amsthm, mathtools, amssymb, color, float}
\usepackage[hyphens]{url}
\usepackage[square,numbers]{natbib}
\usepackage[justification=centering]{caption}

\newcommand{\B}[1]{\textbf{#1}}
\newcommand{\T}{^{\tiny{T}}}
\newcommand{\ds}{\displaystyle}
\newcommand{\Part}[2]{\frac{\partial #1}{\partial #2}}
\newcommand{\PartS}[2]{\frac{\partial^2 #1}{\partial #2^2}}
\newcommand{\N}{\mathcal{N}}

\begin{document}

%\icmlkeywords{Machine Learning, Deep Learning, Neural Networks, Theory of Functional Connections, PDE, ICML}

\title{Deep Theory of Functional Connections: A New Method for Estimating the Solutions of PDEs}
\date{}
\author{Carl Leake\thanks{Ph.D. Student, Aerospace Engineering, Texas A\&M University, College Station, TX. E-mail: \mbox{leakec@tamu.edu}}}

\maketitle

\begin{abstract}
This article presents a new methodology called deep Theory of Functional Connections (TFC) that estimates the solutions of partial differential equations (PDEs) by combining neural networks with TFC. TFC is used to transform PDEs with boundary conditions into unconstrained optimization problems by embedding the boundary conditions into a ``constrained expression.'' In this work, a neural network is chosen as the free function, and used to solve the now unconstrained optimization problem. The loss function is taken as the square of the residual of the PDE. Then, the neural network is trained in an unsupervised manner to solve the unconstrained optimization problem. This methodology has two major differences when compared with popular methods used to estimate the solutions of PDEs. First, this methodology does not need to discretize the domain into a grid, rather, this methodology randomly samples points from the domain during the training phase. Second, after training, this methodology represents a closed form, analytical, differentiable approximation of the solution throughout the entire training domain. In contrast, other popular methods require interpolation if the estimated solution is desired at points that do not lie on the discretized grid. The deep TFC method for estimating the solution of PDEs is demonstrated on four problems with a variety of boundary conditions.
\end{abstract}

\section{Introduction}
Partial differential equations are a powerful mathematical tool that is used to model physical phenomena, and their solutions are used in the design and verification processes of a variety of systems. For example, the solution of PDEs may be used to verify that the loads experienced by a beam throughout its lifetime will not cause the beam to fail. Many methods exist to approximate the solutions of PDEs. The most famous of these methods is the finite element method (FEM) \cite{FEA1,FEA2,FEA3}. FEM has been incredibly successful in approximating the solution to PDEs in a variety of fields including structures, fluids, and acoustics. However, FEM does have some drawbacks. 

FEM discretizes the domain into elements. This works well for low-dimensional cases, but the number of elements grows exponentially with the number of dimensions. Therefore, the discretization becomes prohibitive as the number of dimensions increases. Another issue is that FEM solves the PDE at discrete nodes, but if the solution is needed at locations different than these nodes, an interpolation scheme must be used. 

Reference \cite{ModernPDE} explored the use of neural networks to solve PDEs, and showed that the use of neural networks avoids these problems. Rather than discretizing the entire domain into a number of elements that grows exponentially with the dimension, neural networks can sample points randomly from the domain. Moreover, once the neural network is trained, it is a closed form, analytical, differentiable approximation of the PDE. Therefore, no interpolation scheme is needed when estimating the solution at points that did not appear during training, and further analytical manipulation of the solution can be done with ease. Furthermore, Ref. \cite{ModernPDE} compared the neural network method with FEM on a set of test points that did not appear during training (i.e. points that were not the nodes in the FEM solution), and showed that the solution obtained by the neural network generalized well to points outside of the training data. In fact, the maximum error on the test set of data was never more than the maximum error on the training set of data. In contrast, the FEM had more error on the test set than on the training set. In one case, the test set had approximately 3 orders of magnitude more error than the training set.

Reference \cite{ModernPDE} is not the first article to explore the use of neural networks to solve PDEs. One of the early papers on the topic was Ref. \cite{OrigOdePde}. Although Ref. \cite{ModernPDE} improved what was presented in \cite{OrigOdePde} in almost every way, one item that was missing was exact satisfaction of the boundary constraints. The solution for PDEs presented by Ref. \cite{OrigOdePde} used a method similar to the Coon's Patch \cite{CoonsPatch} to satisfy the boundary constraints exactly.

Exact boundary constraint satisfaction is of interest for a variety of problems, especially when the boundary value information is known to a high degree of precision. This is especially useful for physics informed problems. Take for example, using the heat equation to describe the temperature within a rod. If the temperature on the edges of the rod is know to a high degree of precision, and the heat equation is being used to estimate the temperature inside the rod, then the only information that is known for sure a priori is the boundary conditions. Hence, it is desired to have the boundary conditions of the rod met exactly. Moreover, embedding the boundary conditions means that the neural network only needs to sample points from the interior of the domain, not the domain and the boundary. Thus, a method for embedding boundary conditions for PDEs of arbitrary dimension is desired. Luckily, a framework for embedding constraints into machine learning algorithms has already been invented: the Theory of Functional Connections (TFC) \cite{SVM-TFC}. In Ref. \cite{SVM-TFC}, TFC was used to embed constraints into support vector machines, but left embedding constraints into neural networks to future work. This research shows how to embed constraints into neural networks, and leverages this technique to numerically estimate the solutions of PDEs. 

TFC is a framework that is able to satisfy many types of boundary conditions while maintaining a function that can be freely chosen. This free function can be chosen, for example, to minimize the residual of a differential equation. TFC has already been used to solve differential equations with initial value constraints, boundary value constraints, relative constraints, integral constraints, and infinite constraints \cite{TFC,LDE,NDE,TFCInt}. Recently, the framework was extended to $n$ dimensions \cite{M-TFC} for constraints on the value and flux of $n-1$ dimensional manifolds. This means the framework can now generate constrained expressions that satisfy the boundary constraints of multidimensional PDEs \cite{M-TFC-PDE}. 

\section{Theory of Functional Connections}
The Theory of Functional Connections (TFC) is a mathematical framework designed to turn constrained problems into unconstrained problems. This is accomplished through the use of the \emph{constrained expression}, which is a mathematical expression that represents the family of all possible functions that satisfy the problem's constraints. This technique is especially useful when solving PDEs, as it transforms the PDE into an unconstrained optimization problem. TFC has two major steps: 1) embed the boundary conditions of the problem into the constrained expression 2) solve the now unconstrained optimization problem. The paragraphs that follow will explain these steps in more detail.

The TFC framework is easiest to understand when explained via an example, like a simple harmonic oscillator. Equation \eqref{eq:SHO} gives an example of a simple harmonic oscillator problem.
\begin{equation}\label{eq:SHO}
   my^{\prime\prime}+k y = 0 \quad \text{subject to:}\\
   \begin{cases}
        y(0) = y_0\\
        y^\prime(0) = y^\prime_0
   \end{cases}
\end{equation}

Step one of the TFC framework is to build the constrained expression. The constrained expression consists of two parts; the first part is a function that satisfies the boundary constraints, and the second part projects the free-function, $g(\B{x})$, onto the hyper-surface of functions that are equal to zero at the boundaries. For problems with boundary and derivative constraints, the constrained expression can be written compactly using Eq. (\ref{eq:TFCBound}).
\begin{equation}\label{eq:TFCBound}
\begin{aligned}
   &f(\B{x}) = \mathcal{M}_{i_1,i_2,\dots,i_n}(c(\B{x}))v_{i_1}(x_1)v_{i_2}(x_2)\dots v_{i_n}(x_n) +\\
   &g(\B{x}) - \mathcal{M}_{i_1,i_2,\dots,i_n}(g(\B{x}))v_{i_1}(x_1)v_{i_2}(x_2)\dots v_{i_n}(x_n)
\end{aligned}
\end{equation}
where $\B{x}$ is a vector of the independent variables, $\begin{Bmatrix} x_1 & x_2 & \dots x_n \end{Bmatrix}\T$, $\mathcal{M}$ is an $n$-th order tensor containing the boundary conditions $c(\B{x})$, $v_{i_1}\dots v_{i_n}$ are vectors whose elements are functions of the independent variables, $g(\B{x})$ is the free-function that can be chosen to be anything\footnote{More specifically, the free function $g(\B{x})$ can be chosen to be any function as long as that function is defined at the constraints.}, and $f(\B{x})$ is the constrained expression. The first term in $f(\B{x})$ satisfies the boundary conditions, and the final two terms project the free-function, $g(\B{x})$, onto the space of functions that vanish at the constraints. 

For one dimension with independent variable $x_1$ on the domain $[0,1]$ and for boundary and first-derivative constraints only, $\mathcal{M}$ has the following form,
\begin{equation*}
    \mathcal{M}(c(x_1)) = \begin{bmatrix} c(0) & c_{x_1}(0) & c(1) & c_{x_1}(1) \end{bmatrix}.
\end{equation*}
The values of $\mathcal{M}$ that are unused are eliminated, and the vector, $v$, is created afterwards with the appropriate size. For the one dimensional simple harmonic oscillator, $\mathcal{M}$ and $v$ are given by Eq. (\ref{eq:shoMV}).
\begin{equation}\label{eq:shoMV}
\begin{aligned}
    \mathcal{M}(c(x_1)) &= \begin{bmatrix} c(0) & c_{x_1}(0) \end{bmatrix} \\
    v &= \begin{bmatrix} p^1(x_1) \\ p^2(x_1) \end{bmatrix}
\end{aligned}
\end{equation}
where $p^1({x_1})$ and $p^2(x_1)$ are functions that need to be solved for. Substituting everything into the constrained expression, $f(x_1)$, and simplifying yields Eq. \eqref{eq:shoF}.
\begin{equation}\label{eq:shoF}
   f(x_1) = p^1(x_1)\Big(b(0)-g(0)\Big)+p^2(x_1)\Big(b_{x_1}(0)-g_{x_1}(0)\Big) + g(x_1)
\end{equation}
The functions $p_1(x_1)$ and $p_2(x_1)$ are solved for by setting the constrained expression, $f(x_1)$, equal to the boundary constraints at the boundaries. For the simple harmonic oscillator this yields the set of simultaneous equations given by Eqs. \eqref{eq:shop1} and \eqref{eq:shop2}.
\begin{equation}\label{eq:shop1}
   b(0) = p^1(0)\Big(b(0)-g(0)\Big)+p^2(0)\Big(b_{x_1}(0)-g_{x_1}(0)\Big)+g(0)
\end{equation}
\begin{equation}\label{eq:shop2}
   b_t(0) = p_{x_1}^1(0)\Big(b(0)-g(0)\Big)+p^2_{x_1}(0)\Big(b_{x_1}(0)-g_{x_1}(0)\Big)+g_{x_1}(0) 
\end{equation}
Equation (\ref{eq:shop1}) shows that $p^1(0) = 1$ and $p^2(0)=0$, and Eq. (\ref{eq:shop2}) shows that $p^1_{x_1}(0) = 0$ and $p^2_{x_1}(0)=1$. One solution to this set of simultaneous equations is $p^1(x_1) = 1$ and $p^2(x_1) = x_1$. Substituting these results back into the constrained expression and substituting $b(0)=y_0$ and $b_{x_1}(0)=\dot{y}_0$ yields Eq. (\ref{eq:shof}).
\begin{equation}\label{eq:shof}
    f(x_1) = g(x_1)+y_0-g(0)+t\Big(y^\prime_0-g_{x_1}(0)\Big)
\end{equation}
Notice that for any function $g(x_1)$, the boundary conditions will always be satisfied exactly. Therefore, solving the differential equation has now become an unconstrained optimization problem. This unconstrained optimization problem could be cast in the following way. Let the function to be minimized, $\mathcal{L}$, be equal to the square of the residual of the differential equation,
\begin{equation*}
    \mathcal{L}(x_1)= \Big(mf^{\prime\prime}(x_1)+kf(x_1)\Big)^2.
\end{equation*}
This function is to be minimized by varying the function $g(x_1)$. One way to do this is to make $g(x_1)$ a linear combination of a series of basis functions, and minimize the coefficients that multiply those basis functions via least-squares or some other optimization technique. For an example of this implementation via Chebyshev Orthogonal Polynomials see \cite{LDE}. Also, \cite{NDE} shows how this can be done using non-linear least squares for non-linear ODEs. 

\subsection{$n$-Dimensional Constrained Expressions}
The previous example derived the constrained expression by creating and solving a series of simultaneous algebraic equations. This technique works well for constrained expressions in one dimension; however, a different, more mechanized formalism exists that is useful for deriving these expression in $n$ dimensions for constraints on the value and flux of $n-1$ dimensional manifolds \cite{M-TFC}. This formalism shows how to construct the two portions of the constrained expression shown in Eq. \eqref{eq:TFCBound}: 1) $n$-th order tensor $\mathcal{M}$ 2) the $n$ vectors, $v$. 

Before discussing how to build the $\mathcal{M}$ tensor and $v$ vectors, let's introduce some mathematical notation. Let $k\in[1,n]$ be an index used to denote the $k$-th dimension. Let $^kc_p^d := \frac{\partial^d c(\B{x})}{\partial x_k^d}\bigg|_{x_k=p}$ be the constraint specified by taking the $d$-th derivative of the constraint function, $c(\B{x})$, evaluated at the $x_k=p$ hyperplane. Further, let $^kc_{\B{p}_k}^{\B{d}_k}$ be the vector of $\ell_k$ constraints defined at the $x_k=\B{p}_k$ hyperplanes with derivative orders of $\B{d}_k$, where $\B{p}_k$ and $\B{d}_k \in \mathcal{R}^{\ell_k}$. In addition, let's define a boundary condition operator $^kb_p^d$ that takes the $d$-th derivative with respect to $x_k$ of a function, and then evaluates that function at the $x_k=p$ hyperplane. Mathematically,
\begin{equation*}
    ^kb_p^d[f(\B{x})] = \frac{\partial^d f(\B{x})}{\partial x_k^d}\bigg|_{x_k=p}
\end{equation*}

This mathematical notation will be used to introduce a step-by-step method for building the $\mathcal{M}$ tensor. This step-by-step process will be be illustrated via a 3-dimensional example that has Dirichlet boundary conditions in $x_1$ and initial conditions in $x_2$ and $x_3$ on the domain $x_1,x_2,x_3 \in [0,1]\times[0,1]\times[0,1]$. The $\mathcal{M}$ tensor is constructed using the following three rules.
\begin{enumerate}
    \item The element $\mathcal{M}_{111} = 0$.
    \item The first order sub-tensor of $\mathcal{M}$ specified by keeping one dimension's index free and setting all other dimension's indices to 1 consists of the value $0$ and the boundary conditions for that dimension. Mathematically,
    \begin{equation*}
        \mathcal{M}_{1,\dots,1,i_k,1,\dots,1} = \begin{Bmatrix} 0, ^kc_{\B{p}_k}^{\B{d}_k} \end{Bmatrix}.
    \end{equation*}
    Using the example boundary conditions,
    \begin{align}\label{eq:M1stOrders}
        &\mathcal{M}_{i_111} = \big[0,c(0,x_2,x_3),c(1,x_2,x_3)\big]\T \nonumber\\
        &\mathcal{M}_{1i_21} = \big[0,c(x_1,0,x_3),c_{x_2}(x_1,0,x_3)\big]\T\\
        &\mathcal{M}_{11i_3} = \big[0,c(x_1,x_2,0),c_{x_3}(x_1,x_2,0)\big]\T \nonumber.
    \end{align}
    \item The remaining elements of the $\mathcal{M}$ tensor are those with at least two indices different from one. These elements are the geometric intersection of the boundary condition elements of the first order tensors given in Eq. (\ref{eq:M1stOrders}), plus a sign ($+$ or $-$) that is determined by the number of elements being intersected. Mathematically this can be written as,
    \begin{equation*}
            {\cal M}_{i_1 i_2 \dots i_n} =\, ^1 b^{\B{d}^1_{i_1 - 1}}_{\B{p}^1_{i_1 - 1}} \bigg[\, ^2b^{\B{d}^2_{i_2 - 1}}_{\B{p}^2_{i_2 - 1}} \bigg[ \dots \bigg[\, ^n b^{\B{d}^n_{i_n - 1}}_{\B{p}^n_{i_n - 1}} [c(\B{x})]\bigg] \dots \bigg] \bigg] (-1)^{m+1},
    \end{equation*}
    where $m$ is the number of indices different than one.
    Using the example boundary conditions we give three examples:
    \begin{align*}
       &M_{133} = -c_{x_2x_3}(x_1,0,0)\\
       &M_{221} = -c(0,0,x_3)\\
       &M_{332} = c_{x_2}(1,0,0)
    \end{align*}
\end{enumerate}

A simple procedure also exists for constructing the $v$ vectors. The $v$ vectors have a standard form:
\begin{equation*}
    \B{v}_k = \begin{Bmatrix} 1, & \ds\sum_{i=1}^{\ell_k} \alpha_{i1} \, h_i (x_k), & \ds\sum_{i=1}^{\ell_k} \alpha_{i2} \, h_i (x_k), & \dots, & \ds\sum_{i=1}^{\ell_k} \alpha_{i\ell_k} \, h_i (x_k)\end{Bmatrix},
\end{equation*}
where  $h_i (x_k)$ are $\ell_k$ linearly independent functions. The simplest set of linearly independent functions, and those most often used in the TFC constrained expressions, are monomials, $h_i (x_k) = x_k^{i-1}$. The $\ell_k \times \ell_k$ coefficients, $\alpha_{ij}$, can be computed by matrix~inversion,
\begin{equation*}
    \begin{bmatrix} ^kb^{d_1}_{p_1}[h_1] & ^kb^{d_1}_{p_1}[h_2] & \dots & ^kb^{d_1}_{p_1}[h_{\ell_k}] \\ ^kb^{d_2}_{p_2}[h_1] & ^kb^{d_2}_{p_2}[h_2] & \dots & ^kb^{d_2}_{p_2}[h_{\ell_k}] \\ \vdots & \vdots & \ddots & \vdots \\ ^kb^{d_{\ell_k}}_{p_{\ell_k}}[h_1] & ^kb^{d_{\ell_k}}_{p_{\ell_k}}[h_2] & \dots & ^kb^{d_{\ell_k}}_{p_{\ell_k}}[h_{\ell_k}] \end{bmatrix} \begin{bmatrix} \alpha_{11} & \alpha_{12} & \dots & \alpha_{1\ell_k} \\ \alpha_{21} & \alpha_{22} & \dots & \alpha_{2\ell_k} \\ \vdots & \vdots & \ddots & \vdots \\ \alpha_{\ell_k 1} & \alpha_{\ell_k 2} & \dots & \alpha_{\ell_k\ell_k} \end{bmatrix} = \begin{bmatrix} 1 & 0 & \dots & 0 \\ 0 & 1 & \dots & 0 \\ \vdots & \vdots & \ddots & \vdots \\ 0 & 0 & \dots & 1\end{bmatrix}.
\end{equation*}

Using the example boundary conditions, let's derive the $v_1$ vector using the linearly independent functions $h_1=1$ and $h_2=x_1$.
\begin{equation*}
    \begin{bmatrix} 1 & 0 \\ 1 & 1 \end{bmatrix} \begin{bmatrix} \alpha_{11} & \alpha_{12} \\ \alpha_{21} & \alpha_{22} \end{bmatrix} = \begin{bmatrix} 1 & 0\\ 0 & 1\end{bmatrix} \implies \begin{bmatrix} \alpha_{11} & \alpha_{12} \\ \alpha_{21} & \alpha_{22} \end{bmatrix} = \begin{bmatrix} 1 & 0 \\ -1 & 1 \end{bmatrix}
\end{equation*}
\begin{equation*}
    v_1 = \begin{Bmatrix}1 & 1-x_1 & x_1\end{Bmatrix}.
\end{equation*}

\noindent For more examples and a proof that these procedures for generating the $\mathcal{M}$ tensor and $v$ vectors form a constrained expression see Ref. \cite{M-TFC}.

\subsection{Two-Dimensional Example}

This subsection will give an in depth example for a two-dimensional TFC case. The example is originally from \cite{OrigOdePde} problem 5, and is one of the PDE problems analyzed in this article. The problem is shown in Eq. (\ref{eq:2dEx}).
\begin{equation}\label{eq:2dEx}
\begin{aligned}
   &\nabla^2 z(x,y) = e^{-x}(x-2+y^3+6y)\quad  \text{subject to:}\\
   &\begin{cases}
   z(x,0) = c(x,0)= xe^{-x}\\
   z(0,y) = c(0,y) = y^3\\
   z(x,1) = c(x,1) = e^{-x}(x+1)\\
   z(1,y) = c(1,y) = (1+y^3)e^{-1}
   \end{cases}\\
   &\text{where} \quad (x,y) \in [0,1]\times[0,1]
\end{aligned}
\end{equation}

Following the step-by-step procedure given in the previous section we will construct the $\mathcal{M}$ tensor:
\begin{enumerate}
    \item The first element is $\mathcal{M}_{11} = 0$.
    \item The first order sub-tensors of $\mathcal{M}$ are: 
    \begin{align*}
        \mathcal{M}_{i_11} &= \begin{Bmatrix} 0 & c(0,y) & c(1,y)\end{Bmatrix}\\
        \mathcal{M}_{1i_2} &= \begin{Bmatrix} 0 & c(x,0)& c(x,1)\end{Bmatrix}
    \end{align*}
    \item The remaining elements of $\mathcal{M}$ are the geometric intersection of elements from the first order sub-tensors. 
    \begin{align*}
        \mathcal{M}_{22} &= -c(0,0) &&\mathcal{M}_{23} = -c(1,0) \\
        \mathcal{M}_{32} &= -c(0,1) &&\mathcal{M}_{33} = -c(1,1)
    \end{align*}
\end{enumerate}
Hence, the $\mathcal{M}$ tensor is,
\begin{equation*}
    \mathcal{M}_{i_1i_2} = \begin{bmatrix} 0 & c(0,y) & c(1,y) \\ c(x,0) & -c(0,0) & -c(1,0) \\ c(x,1) & -c(0,1) & -c(1,1)\end{bmatrix} = \begin{bmatrix} 0 & y^3 & (1+y^3)e^{-1} \\ xe^{-x} & 0 & -e^{-1} \\ e^{-x}(x+1) & -1 & -2e^{-1} \end{bmatrix} 
\end{equation*}

Following the step-by-step procedure given in the previous section we will construct the $v$ vectors. For $v_{i_1}$ let's choose the linearly independent functions $h_1=1$ and $h_2=x$.
\begin{equation*}
    \begin{bmatrix} 1 & 0 \\ 1 & 1 \end{bmatrix} \begin{bmatrix} \alpha_{11} & \alpha_{12} \\ \alpha_{21} & \alpha_{22} \end{bmatrix} = \begin{bmatrix} 1 & 0\\ 0 & 1\end{bmatrix} \implies \begin{bmatrix} \alpha_{11} & \alpha_{12} \\ \alpha_{21} & \alpha_{22} \end{bmatrix} = \begin{bmatrix} 1 & 0 \\ -1 & 1 \end{bmatrix}
\end{equation*}
\begin{equation*}
    v_{i_1} = \begin{Bmatrix}1 & 1-x & x\end{Bmatrix}.
\end{equation*}
For $v_{i_2}$ let's choose the linearly indpendent functions $h_1=1$ and $h_2=y$.
\begin{equation*}
    \begin{bmatrix} 1 & 0 \\ 1 & 1 \end{bmatrix} \begin{bmatrix} \alpha_{11} & \alpha_{12} \\ \alpha_{21} & \alpha_{22} \end{bmatrix} = \begin{bmatrix} 1 & 0\\ 0 & 1\end{bmatrix} \implies \begin{bmatrix} \alpha_{11} & \alpha_{12} \\ \alpha_{21} & \alpha_{22} \end{bmatrix} = \begin{bmatrix} 1 & 0 \\ -1 & 1 \end{bmatrix}
\end{equation*}
\begin{equation*}
    v_{i_2} = \begin{Bmatrix}1 & 1-y & y\end{Bmatrix}.
\end{equation*}

Now, we use the constrained expression form given in Eq. \eqref{eq:TFCBound} to finish building the constrained expression.
\begin{align}\label{eq:fExSoln}
    f(x,y) = \ &g(x,y)+\frac{xy(y^2-1)}{e}+e^{-x}(x+y)+(1-x)\bigg(g(0,0)+y\Big(g(0,1)+y^2-g(0,0)-1\Big)\bigg)+\\ &(x-1)g(0,y)+x\Big(yg(1,1)+(1-y)g(1,0)\Big)-xg(1,y)+(y-1)g(x,0)-yg(x,1)\nonumber 
\end{align}

Notice, that Eq. (\ref{eq:fExSoln}) will always satisfy the boundary conditions of the problem regardless of the value of $g(x,y)$. Thus, the problem has been transformed into an unconstrained optimization problem where the cost function, $\mathcal{L}$, is the square of the residual of the PDE, 
\begin{equation*}
    \mathcal{L}(x,y) = \Big(\nabla^2f(x,y)-e^{-x}(x-2+y^3+6y)\Big)^2.
\end{equation*}
For one-dimensional ODEs, the minimization of the cost function was accomplished by making $g$ the summation of orthogonal polynomials, and performing least-squares or some other optimization technique to find the coefficients that multiply those orthogonal polynomials. For two dimensions, one could make $g(x,y)$ the product of two sums of these orthogonal polynomials, calculate all of the cross-terms, and then solve for the coefficients that multiply all terms and cross-terms using least-squares or non-linear least-squares. However, this will become prohibitive as the dimension increases. An alternative solution, and the one explored in this article, is to make the free function, $g(x,y)$, a neural network.

\section{PDE Solution Methodology}
As with the previous section, the easiest way to describe the methodology is with an example. The example used throughout this section will be the PDE given in Eq. (\ref{eq:2dEx}).

As mentioned in the previous section, deep TFC approximates solutions to PDEs by finding the constrained expression for the PDE and choosing a neural network as the free function. For all of the problems analyzed in this article, a simple, fully connected neural network was used. These networks consist of non-linear activation functions composed with affine transformations of the form $\mathcal{A} = W\cdot x + b$, where $W$ are neuron weights, $b$ are neuron biases, and $x$ is a vector of inputs from the previous layer (or the inputs to the neural network if it is the first layer). The weights and biases of the entire neural network make up the tunable parameters, $\theta$. In this article, neural networks will be represented with the symbol $\mathcal{N}$. For example, a neural network with inputs $x$ and $y$ would be given as $\mathcal{N}(x,y;\theta)$. Thus, the constrained expression, given originally in Eq. (\ref{eq:fExSoln}), now has the form given in Eq. (\ref{eq:conPDE}).
\begin{equation}\label{eq:conPDE}
\begin{aligned}
    f(x,y) = \ &\mathcal{N}(x,y;\theta)+\frac{xy(y^2-1)}{e}+e^{-x}(x+y)+\\
    &(1-x)\bigg(\mathcal{N}(0,0;\theta)+y\Big(\mathcal{N}(0,1;\theta)+y^2-\mathcal{N}(0,0;\theta)-1\Big)\bigg)+(x-1)\mathcal{N}(0,y;\theta)+\\
    &x\Big(y\mathcal{N}(1,1;\theta)+(1-y)\mathcal{N}(1,0;\theta)\Big)-x\mathcal{N}(1,y;\theta)+(y-1)\mathcal{N}(x,0;\theta)-y\mathcal{N}(x,1;\theta)
\end{aligned}
\end{equation}

In order to estimate the solution to the PDE, the parameters of the neural network have to be optimized to minimize the loss function, which is taken to be the square of the residual of the PDE. For this example,
\begin{align*}
    \mathcal{L}_i(x_i,y_i) &= \Big(\nabla^2f(x_i,y_i)-e^{-x_i}(x_i-2+y_i^3+6y_i)\Big)^2,\\
    \mathcal{L} &= \sum_i^N \mathcal{L}_i.
\end{align*}
The attentive reader will notice that training the neural network will require, for this example, taking two second order partial derivatives of $f(x,y)$ to calculate $\mathcal{L}_i$, and then taking gradients of $\mathcal{L}$ with respect to the neural network parameters, $\theta$, in order to train the neural network.

To take these higher order derivatives, TensorFlow's\texttrademark\ gradients function was used \cite{tensorflowWhitepaper}. This function uses automatic differentiation \cite{AD} to compute these derivatives. However, one must be conscientious when using the gradients function to ensure they get the correct gradients. 

When taking the gradient of a vector, $y_j$, with respect to another vector, $x_i$, TensorFlow\texttrademark\ computes,
\begin{equation*}
    z_i = \frac{\partial \big(\sum_{j=1}^N y_j\big)}{\partial x_i},
\end{equation*}
where $z_i$ is a vector of the same size as $x_i$. The only place where this may be an issue in the example used in this section is when computing $\nabla^2 f_i$. The desired output of this calculation is the following vector,
\begin{equation*}
    z_i = \bigg\{\frac{\partial^2 f_1}{\partial x_1^2}+\frac{\partial^2 f_1}{\partial y_1^2},\dots,\frac{\partial^2 f_N}{\partial x_N^2}+\frac{\partial^2 f_N}{\partial y_N^2}\bigg\},
\end{equation*}
where $z_i$ has the same size as $f_i$ and $(x_i,y_i)$ is the pair used to generate $f_i$. TensorFlow's\texttrademark\ gradients function will compute the following vector,
\begin{align*}
    z_i = \bigg\{&\frac{\partial^2 \big(\sum_{j=1}^N f_j \big)}{\partial x_1^2}+\frac{\partial^2 \big(\sum_{j=1}^N f_j \big)}{\partial y_1^2},\dots,\\
    &\frac{\partial^2 \big(\sum_{j=1}^N f_j \big)}{\partial x_N^2}+\frac{\partial^2 \big(\sum_{j=1}^N f_j \big)}{\partial y_N^2}\bigg\}.
\end{align*}
However, because $f_i$ only depends on $(x_i,y_i)$ and the derivative operator commutes with the sum operator, TensorFlow's\texttrademark\ gradients function will compute the desired vector. Moreover, the size of the output vector will be correct because the input vectors, $x_i$ and $y_i$, have the same size as $f_i$.

\subsection{Training the Neural Network}
Three methods were tried when optimizing the parameters of the neural networks:
\begin{enumerate}
    \item Adam optimizer \cite{Adam}
    \item Broyden-–Fletcher–-Goldfarb-–Shanno (BFGS) optimizaion algorithm \cite{OrigTrainMethod}
    \item Hybrid method
\end{enumerate}
The first method, Adam, is a variant of stochastic gradient descent (SGD) that combines the advantages of two other popular SGD variants: AdaGrad \cite{AdaGrad} and RMSProp \cite{RMSProp}.

The second method, BFGS, is a quasi-Newton method designed for solving unconstrained, non-linear optimization problems. This method was chosen based on its performance when optimizing neural network parameters to estimate PDE solutions in Ref. \cite{OrigOdePde}. The hybrid method uses both methods in series. For all four problems shown in this article, the BFGS optimizer gave the best results. 

\section{Results}
This section compares the estimated solution found using deep TFC with with the analytical solution. Four PDE problems are analyzed. The first is the example PDE given in Eq. (\ref{eq:2dEx}), and the second is the wave equation. The third and fourth PDEs are simple solutions to the incompressible Navier-Stokes equations.

\subsection{Problem 1}
The first problem analyzed was the PDE given by Eq. (\ref{eq:2dEx}), copied below for the readers' convenience.
\begin{equation*}
\begin{aligned}
   &\nabla^2 z(x,y) = e^{-x}(x-2+y^3+6y)\quad  \text{subject to:}\\
   &\begin{cases}
   z(x,0) = xe^{-x}\\
   z(0,y) = y^3\\
   z(x,1) = e^{-x}(x+1)\\
   z(1,y) = (1+y^3)e^{-1}
   \end{cases}\\
   &\text{where} \quad (x,y) \in [0,1]\times[0,1]
\end{aligned}
\end{equation*}
The neural network used to estimate the solution to this PDE was a fully connected neural network with five hidden layers and 30 neurons per layer. The non-linear activation function used was the hyperbolic tangent. Other neural network sizes and non-linear activation functions were tried, but this size and activation function combination did the best. The biases of the neural network were all initialized as zeros, and the weights were initialized using TensorFlow's\texttrademark\ implementation of the Xavier initialization with uniform random initialization \cite{Init}. Training pairs, $(x,y)$, were created for this problem by sampling $x$ and $y$ as independent and identically distributed (IID) random variables on the uniform distribution on $[0,1]$. The network was trained using the BFGS method on a batch size of 10,000 training pairs. 

\begin{figure}[!h]
    \centering
    \includegraphics[trim={2.5cm 1.0cm 0.75cm 2.0cm},clip,width=0.6\linewidth]{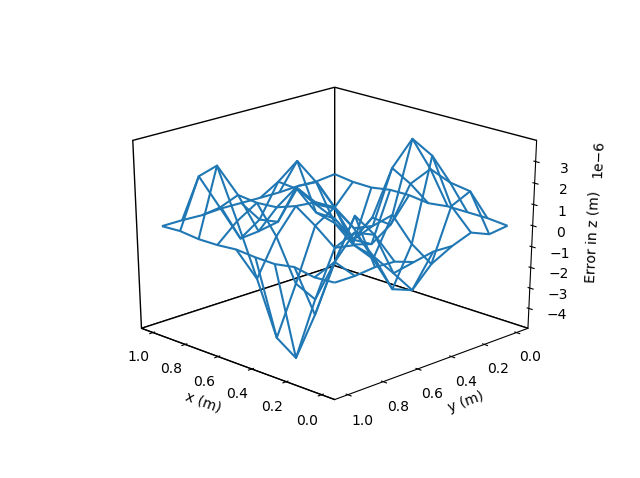}
    \caption{Problem 1 Solution Error}
    \label{fig:MyP5}
\end{figure}

Figure \ref{fig:MyP5} shows the difference between the analytical solution and the estimated solution using deep TFC on a grid of 100 evenly distributed points (10 per independent variable). This grid represents the test set. 

The maximum error on the test set was $4.806\times 10^{-6}$ meters and the average error was $8.872\times 10^{-7}$ meters. The maximum error is relatively low, five orders of magnitude lower than the solution values, which are on the order of $10^{-1}$ meters. However, the maximum solution error using deep TFC is not as low as the maximum solution error obtained in \cite{OrigOdePde}, which was on the order of $10^{-7}$ meters. Although there are many possible explanations for this discrepancy, for example, different implementations of the optimizer or different weight initialization functions, the author was concerned that it may be the assumed solution form, $f(x,y)$, that was causing the increase in estimated solution error. The solution form created using deep TFC is more complex both in the number of terms and the number of times the neural network appears within the assumed solution form.

To investigate, a comparison was made between the solution form posed in \cite{OrigOdePde} and the solution form posed in this article, while keeping all other variables constant. For this comparison, the neural network architecture consisted of one hidden layer with 10 neurons, and the non-linear activation function used was the hyperbolic tangent function\footnote{The sigmoid function was also tried, which is what is used in \cite{OrigOdePde}; however, the average error and maximum error for both solution forms, deep TFC and \cite{OrigOdePde}, was almost doubled compared to using the hyperbolic tangent function.}. The neural network was trained using the BFGS optimizer. The training pairs were created by randomly sampling 10,000 point pairs $(x,y)$ where $x$ and $y$ are IID random variables on the uniform distribution $[0,1]$. The test set was a grid of 100 evenly distributed points (10 per independent variable).

\begin{figure}[h]
    \centering
    \includegraphics[trim={2.5cm 1.0cm 0.75cm 2.0cm},clip,width=0.6\linewidth]{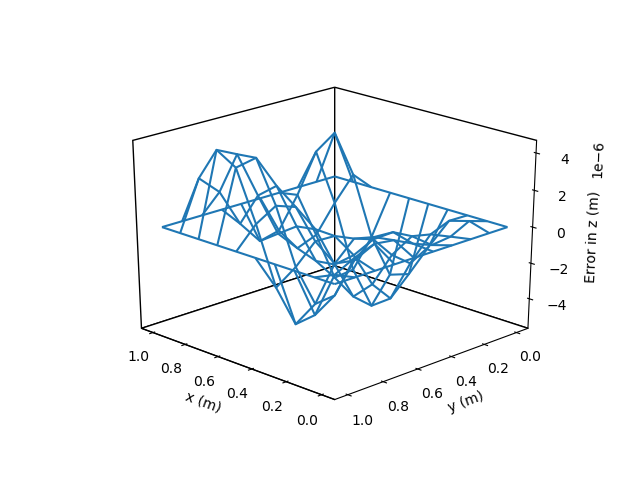}
    \caption{Problem 1 Solution Error Using \cite{OrigOdePde} Solution Form}
    \label{fig:P5OldMethod}
\end{figure}

Figure \ref{fig:P5OldMethod} was created using the solution form posed in \cite{OrigOdePde}. The maximum error on the test set was $5.481\times 10^{-6}$ meters and the average error on the test set was $1.131\times 10^{-6}$ meters.

\begin{figure}[h]
    \centering
    \includegraphics[trim={2.5cm 1.0cm 0.75cm 2.0cm},clip,width=0.6\linewidth]{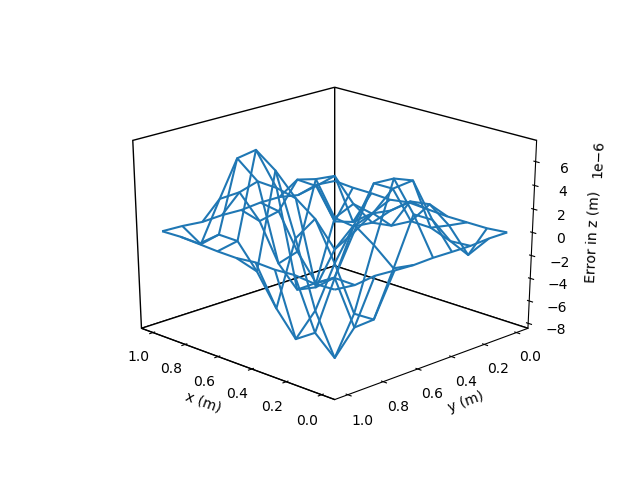}
    \caption{Problem 1 Solution Error Using Deep TFC Solution Form}
    \label{fig:P5TFC}
\end{figure}

Figure \ref{fig:P5TFC} was created using the deep TFC solution form. The maximum error on the test set was $8.124\times 10^{-6}$ meters and the average error on the test set was $1.696\times 10^{-6}$ meters.

Comparing Figs. \ref{fig:P5OldMethod} and \ref{fig:P5TFC} shows that the solution form from \cite{OrigOdePde} does slightly better in terms of average error and maximum error for this problem, but the difference is not significant. However, neither of the tests conducted here is able to reduce the error to the level reported in \cite{OrigOdePde}. Thus, there must be some other factor, initialization, optimizer method, etc., that is causing this difference. One important take-away from this comparison is that despite the more complex solution form created using deep TFC, which can easily be applied to higher dimensions, the final solution accuracy is very similar to the simpler solution form that was designed for lower-dimensional problems. 

\subsection{Problem 2}
The second problem analyzed was the wave equation, shown in Eq. (\ref{eq:Prob2}).
\begin{equation}\label{eq:Prob2}
\begin{aligned}
    &\frac{\partial^2 u}{\partial t^2}(x,t) = c^2 \frac{\partial^2 u}{\partial x^2}(x,t) \quad \text{subject to:}\\
    &\begin{cases}
    &u(0,t)=0\\
    &u(1,t)=0\\
    &u(x,0)=x(1-x)\\
    &u_t(x,0)=0
    \end{cases}\\
    &\text{where} \quad (x,t) \in [0,1] \times [0,1]
\end{aligned}
\end{equation}
where the constant, $c$, was chosen to be 1. The constrained expression for this problem is shown in Eq. (\ref{eq:TFCWave}).
\begin{equation}\label{eq:TFCWave}
\begin{aligned}
   f(x,t) =\ & (1-x)\Big[\N(0,0)-\N(0,t)\Big]+x\Big[\N(1,0)-\N(1,t)\Big]-\N(x,0)+x(1-x)+\N(x,t)+\\
   &t\Big[(1-x)\N_t(0,0)+x\N_t(1,0)-\N_t(x,0)\Big]
\end{aligned}
\end{equation}
The neural network used to estimate the solution to this PDE was a fully connected neural network with two hidden layers and 30 neurons per layer. The non-linear activation function used was the hyperbolic tangent. The biases and weights were initialized using the same method as problem 1. The training points, $(x,t)$, were created by sampling $x$ and $t$ IID from the uniform distribution on $[0,1]$. The network was trained using the BFGS method on a batch size of 10,000 training pairs.

\begin{figure}[H]
    \centering
    \includegraphics[trim={2.5cm 1.0cm 0.75cm 2.0cm},clip,width=0.6\linewidth]{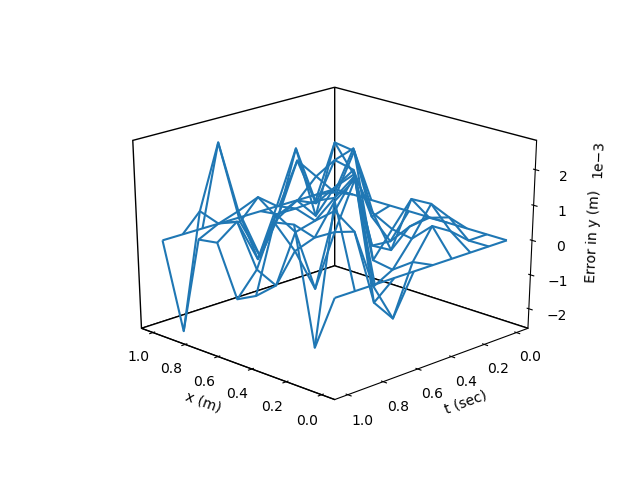}
    \caption{Problem 2 Solution Error}
    \label{fig:P2}
\end{figure}

Figure \ref{fig:P2} shows the difference between the analytical solution and the estimated solution using deep TFC on a grid of 100 evenly distributed points (10 per independent variable). This grid represents the test set. 

The maximum error on the test set was $2.701\times 10^{-3}$ meters and the average error on the test set was $6.978\times 10^{-4}$ meters. The error of this solution is larger than in the problem one, while the solution values are on the same order of magnitude, $10^{-1}$ meters, as in problem one. The larger relative error in problem two is most likely due to the more oscillatory nature of the solution. 

\subsection{Problem 3}
The third problem analyzed was a known solution to the incompressible Navier-Stokes equations, called Poiseuille flow. The problem solves the flow velocity in a two-dimensional pipe in steady-state with a constant pressure gradient applied in the longitudinal axis. Equation (\ref{eq:NS1}) shows the associated equations and boundary conditions.
\begin{equation}\label{eq:NS1}
\begin{aligned}
    &\Part{u}{x}+\Part{v}{y} = 0\\
    &\rho \bigg(\Part{u}{t}+u\Part{u}{x}+v\Part{u}{y}\bigg) = -\Part{P}{x}+\mu \bigg(\PartS{u}{x}+\PartS{u}{y}\bigg)\\
    &\rho \bigg(\Part{v}{t}+u\Part{v}{x}+v\Part{v}{y}\bigg) = \mu \bigg(\PartS{v}{x}+\PartS{v}{y}\bigg)\\
    &\text{subject to:}\\
    &\begin{cases}
    &u(0,y,t) = u(L,y,t) = u(x,y,0) = \frac{1}{2\mu}\Part{P}{x}\bigg(y^2-\Big(\frac{H}{2}\Big)^2\bigg)\\
    &u(x,\frac{H}{2},t) = u(x,-\frac{H}{2},t) = 0\\
    &v(0,y,t)=v(L,y,t)=v(x,y,0)=0\\
    &v(0,\frac{H}{2},t)=v(0,-\frac{H}{2},t)=0
    \end{cases}
\end{aligned}
\end{equation}
where $u$ and $v$ are velocities in the $x$ and $y$ directions respectively, $H$ is the height of the channel, $P$ is the pressure, $\rho$ is the density, and $\mu$ is the viscosity. For this problem, the values $H=1$ m, $\rho=1$ kg/m$^3$, $\mu=1$ Pa$\cdot$s, and $\Part{P}{x}=-5$ N/m$^3$ were chosen. The constrained expressions for the $u$-velocity, $f^u(x,y,t)$, and $v$-velocity, $f^v(x,y,t)$, are shown in Eq. (\ref{eq:TFCNS}).
\begin{equation}\label{eq:TFCNS}
\begin{aligned}
   f^u(x,y,t) =\ & \N(x,y,t;\theta)-\N(x,y,0;\theta)+\frac{L-x}{L}\bigg(\N(0,y,0;\theta)-\N(0,y,t;\theta)\bigg)\\
   &+\frac{x}{L}\bigg(\N(L,y,0;\theta)-\N(L,y,t;\theta)\bigg)+\frac{P\Big(4y^2-H^2\Big)}{8\mu}+\\
   &\frac{1}{2HL}\Bigg( (2y-H)\bigg((L-x)\N\Big(0,-\frac{H}{2},0;\theta\Big)+x\N\Big(L,-\frac{H}{2},0;\theta\Big)-L\N\Big(x,-\frac{H}{2},0;\theta\Big)\\
   &-(L-x)\N\Big(0,-\frac{H}{2},t;\theta\Big)+L\N\Big(x,-\frac{H}{2},t;\theta\Big)-x\N\Big(L,-\frac{H}{2},t;\theta\Big)\bigg)\\
   &-(H+2y)\bigg((L-x)\N\Big(0,\frac{H}{2},0;\theta\Big)-L\N\Big(x,\frac{H}{2},0;\theta\Big)+x\N\Big(L,\frac{H}{2},0;\theta\Big)\\
   &-(L-x)\N\Big(0,\frac{H}{2},t;\theta\Big)-x\N\Big(L,\frac{H}{2},t;\theta\Big)+L\N\Big(x,\frac{H}{2},t;\theta\Big)\bigg)\Bigg)\\
   f^v(x,y,t) =\ & \N(x,y,t;\theta)-\N(x,y,0;\theta)+\frac{L-x}{L}\bigg(\N(0,y,0;\theta)-\N(0,y,t;\theta)\bigg)\\
   &+\frac{x}{L}\bigg(\N(L,y,0;\theta)-\N(L,y,t;\theta)\bigg)+\frac{1}{2HL}\Bigg( (2y-H)\bigg((L-x)\N\Big(0,-\frac{H}{2},0;\theta\Big)\\
   &+x\N\Big(L,-\frac{H}{2},0;\theta\Big)-L\N\Big(x,-\frac{H}{2},0;\theta\Big)-(L-x)\N\Big(0,-\frac{H}{2},t;\theta\Big)+L\N\Big(x,-\frac{H}{2},t;\theta\Big)\\
   &-x\N\Big(L,-\frac{H}{2},t;\theta\Big)\bigg)-(H+2y)\bigg((L-x)\N\Big(0,\frac{H}{2},0;\theta\Big)-L\N\Big(x,\frac{H}{2},0;\theta\Big)\\
   &+x\N\Big(L,\frac{H}{2},0;\theta\Big)-(L-x)\N\Big(0,\frac{H}{2},t;\theta\Big)-x\N\Big(L,\frac{H}{2},t;\theta\Big)+L\N\Big(x,\frac{H}{2},t;\theta\Big)\bigg)\Bigg)\\
\end{aligned}
\end{equation}
The neural network used to estimate the solution to this PDE was a fully connected neural network with four hidden layers and 30 neurons per layer. The non-linear activation function used was the sigmoid. The biases and weights were initialized using the same method as problem 1. The training points, $(x,y,t)$, were created by sampling $x$, $y$, and $t$ IID from the a uniform distribution that spanned the range of the associated independent variable. For $x$, the range was $[0,1]$. For $y$, the range was $[-\frac{H}{2},\frac{H}{2}]$, and for $t$, the range was $[0,1]$. The network was trained using the BFGS method on a batch size of 1,000 training pairs. The loss function used was the sum of the squares of the residuals of the three PDEs in Eq. (\ref{eq:NS1}).

For one particular run\footnote{All errors are given in meters per second.}, the maximum error in the $u$-velocity was $3.308\times 10^{-7}$, the average error in the $u$-velocity was $9.998\times 10^{-8}$, the maximum error in the $v$-velocity was $5.575 \times 10^{-7}$, and the average error in the $v$-velocity was $1.542\times 10^{-7}$. The maximum error and average error for this problem are much lower than in problems 1 and 2. However, the constrained expression for this problem essentially encodes the solution, because the initial flow condition at time zero is the same as the flow condition throughout the spatial domain at any time. Thus, if the neural network outputs a value of zero for all inputs, the problem will be solved exactly. Although the neural network does output a very small value for all inputs, it is interesting to note that none of the layers have weights or biases that are at or near zero.

\subsection{Problem 4}
The fourth problem is another solution to the Navier-Stokes equations, and is very similar to the third. The only difference is that in this case, the fluid is not in steady state, it starts from rest. Equation \eqref{eq:NS2} shows the associated equations and boundary conditions.
\begin{equation}\label{eq:NS2}
\begin{aligned}
    &\Part{u}{x}+\Part{v}{y} = 0\\
    &\rho \bigg(\Part{u}{t}+u\Part{u}{x}+v\Part{u}{y}\bigg) = -\Part{P}{x}+\mu \bigg(\PartS{u}{x}+\PartS{u}{y}\bigg)\\
    &\rho \bigg(\Part{v}{t}+u\Part{v}{x}+v\Part{v}{y}\bigg) = \mu \bigg(\PartS{v}{x}+\PartS{v}{y}\bigg)\\
    &\text{subject to:}\\
    &\begin{cases}
    &u(0,y,t) = \Part{u}{x}(0,y,t) = u(x,y,0) = 0\\
    &u(x,\frac{H}{2},t) = u(x,-\frac{H}{2},t) = 0\\
    &v(0,y,t)=\Part{v}{x}(0,y,t)=v(x,y,0)=0\\
    &v(0,\frac{H}{2},t)=v(0,-\frac{H}{2},t)=0
    \end{cases}
\end{aligned}
\end{equation}
This problem was created to avoid encoding the solution to the problem into the constrained expression, as was the case in the previous problem. 
The constrained expressions for the $u$-velocity, $f^u(x,y,t)$, and $v$-velocity, $f^v(x,y,t)$, are shown in Eq. (\ref{eq:TFCNS2}).
\begin{equation}\label{eq:TFCNS2}
\begin{aligned}
   f^u(x,y,t) =\ & \N(x,y,t;\theta)-\N(x,y,0;\theta)+\N(0,y,0;\theta)-\N(0,y,t;\theta)+x\N_x(0,y,0;\theta)-x\N_x(0,y,t;\theta)\\
   &\frac{1}{2H}\Bigg( (2y-H)\bigg(\N\Big(0,-\frac{H}{2},0;\theta\Big)-\N\Big(x,-\frac{H}{2},0;\theta\Big)+x\N_x\Big(0,-\frac{H}{2},0;\theta\Big)-\N\Big(0,-\frac{H}{2},t;\theta\Big)\\
   &+\N\Big(x,-\frac{H}{2},t;\theta\Big)-x\N_x\Big(0,-\frac{H}{2},t;\theta\Big)\bigg)-(H+2y)\bigg(\N\Big(0,\frac{H}{2},0;\theta\Big)-\N\Big(x,\frac{H}{2},0;\theta\Big)\\
   &+x\N_x\Big(0,\frac{H}{2},0;\theta\Big)-\N\Big(0,\frac{H}{2},t;\theta\Big)+\N\Big(x,\frac{H}{2},t;\theta\Big)-x\N_x\Big(0,\frac{H}{2},t;\theta\Big)\bigg)\Bigg)\\
   f^v(x,y,t) =\ & \N(x,y,t;\theta)-\N(x,y,0;\theta)+\N(0,y,0;\theta)-\N(0,y,t;\theta)+x\N_x(0,y,0;\theta)-x\N_x(0,y,t;\theta)\\
   &\frac{1}{2H}\Bigg( (2y-H)\bigg(\N\Big(0,-\frac{H}{2},0;\theta\Big)-\N\Big(x,-\frac{H}{2},0;\theta\Big)+x\N_x\Big(0,-\frac{H}{2},0;\theta\Big)-\N\Big(0,-\frac{H}{2},t;\theta\Big)\\
   &+\N\Big(x,-\frac{H}{2},t;\theta\Big)-x\N_x\Big(0,-\frac{H}{2},t;\theta\Big)\bigg)-(H+2y)\bigg(\N\Big(0,\frac{H}{2},0;\theta\Big)-\N\Big(x,\frac{H}{2},0;\theta\Big)\\
   &+x\N_x\Big(0,\frac{H}{2},0;\theta\Big)-\N\Big(0,\frac{H}{2},t;\theta\Big)+\N\Big(x,\frac{H}{2},t;\theta\Big)-x\N_x\Big(0,\frac{H}{2},t;\theta\Big)\bigg)\Bigg)\\
\end{aligned}
\end{equation}

\begin{figure*}
\begin{minipage}{0.33\linewidth}
    \centering
    \includegraphics[trim={0.0cm 0.0cm 1.25cm 1.25cm},clip,width=\linewidth]{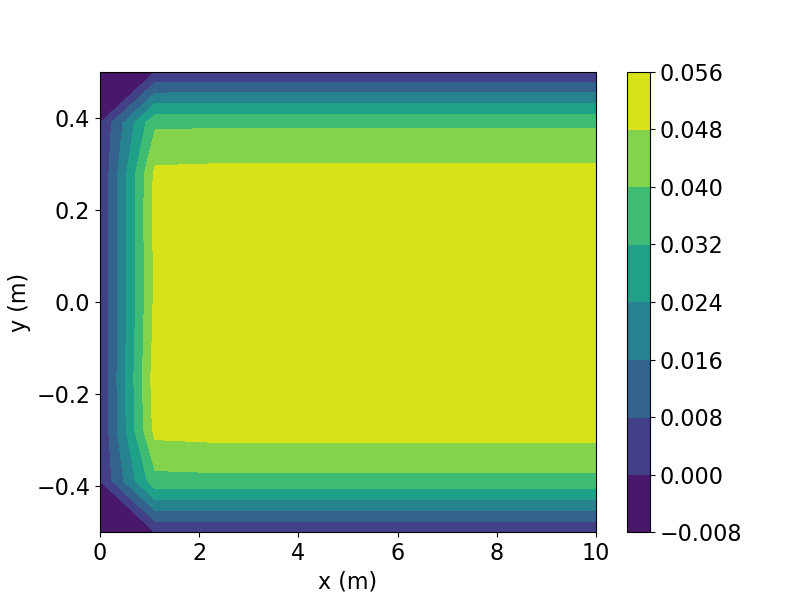}
    \parbox{0.70\linewidth}{\caption{U-Velocity in Meters Per Second at 0.01 Seconds}
    \label{fig:Prob4_0.01}}
\end{minipage}%
\begin{minipage}{0.33\linewidth}
    \centering
    \includegraphics[trim={0.0cm 0.0cm 1.25cm 1.25cm},clip,width=\linewidth]{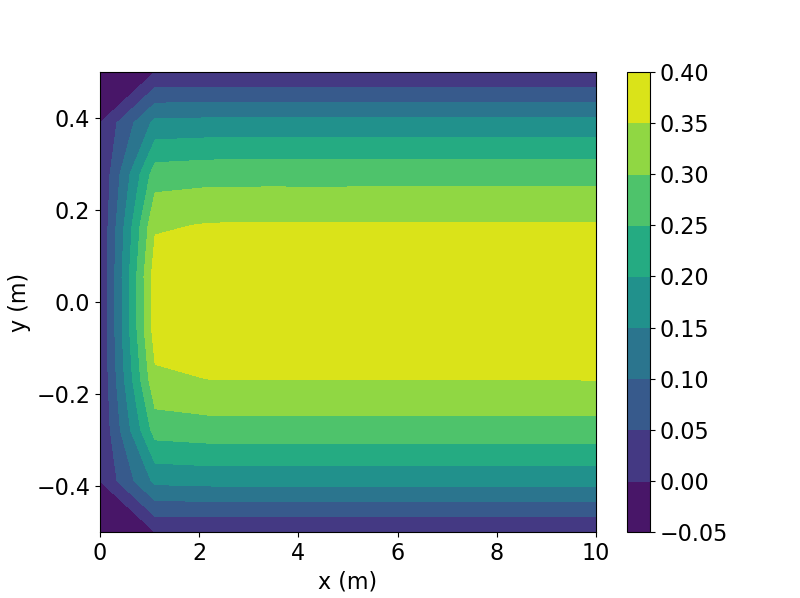}
    \parbox{0.70\linewidth}{\caption{U-Velocity in Meters Per Second at 0.1 Seconds}
    \label{fig:Prob4_0.1}}
\end{minipage}%
\begin{minipage}{0.33\linewidth}
    \centering
    \includegraphics[trim={0.0cm 0.0cm 1.25cm 1.25cm},clip,width=\linewidth]{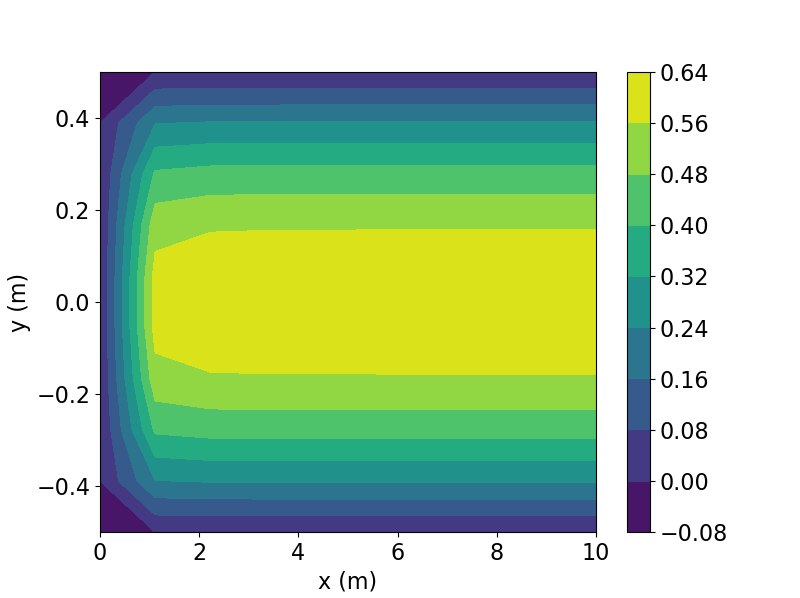}
    \parbox{0.70\linewidth}{\caption{U-Velocity in Meters Per Second at 3.0 Seconds}
    \label{fig:Prob4_1.0}}
\end{minipage}%
\end{figure*}

The neural network used in this problem is exactly the same as the neural network used in problem three. Problem 4 used 2,000 training points that were selected the same way as in problem three, except the new ranges for the independent variables were $[0,10]$ for $x$, $[0,3]$ for $t$, and $[-\frac{H}{2},\frac{H}{2}]$ for $y$.  

Figures \ref{fig:Prob4_0.01} through \ref{fig:Prob4_1.0} show the $u$-velocity of the fluid throughout the domain at three different times. Qualitatively, the solution should look as follows. The solution should be symmetric about the line $y=0$, and the solution should develop spatially and temporally such that after a period of time and sufficiently far from the inlet, $x=0$, the $u$-velocity will be equal or very close to the steady state $u$-velocity of problem 3. Qualitatively, the $u$-velocity field looks correct at all points. Quantitatively, the $u$-velocity at $x=10$ from Fig. \ref{fig:Prob4_1.0} was compared with the known steady state $u$-velocity, and had a maximum error of $5.530\times 10^{-4}$ meters per second and an average error of $2.742\times 10^{-4}$ meters per second.

\section{Conclusions}
This article demonstrated how to combine neural networks with the Theory of Functional Connections into a new methodology, called deep TFC, that was used to estimate the solutions of PDEs. Results on four problems were presented that display how accurately relatively simple neural networks can approximate the solutions to some well known PDEs. The difficulty of the PDEs in these problems ranged from linear, two-dimensional PDEs to coupled, non-linear, three-dimensional PDEs. Moreover, while the focus of this article was on PDEs, the capability to embed constraints into neural networks is useful for many other applications as well.

Future work should investigate the performance of different neural network architectures on the estimated solution error. For example, Ref. \cite{ModernPDE} suggests a neural network architecture where the hidden layers contain element-wise multiplications and sums of sub-layers. The sub-layers are more standard neural network layers like the fully connected layers used in the neural networks of this article.

Another topic for investigation is reducing the estimated solution error by sampling the training points based on the loss function values for the training points of the previous iteration. For example, one could create batches where half of the new batch consists of half of the points in the previous batch that had the largest loss function value and the other half are randomly sampled from the domain. This should consistently give training points that are in portions of the domain where the estimated solution is farthest from the real solution.

\section*{Acknowledgements}
Carl Leake would like to acknowledge the support of the NASA Space Technology Research Fellowship (NSTRF) 2019: grant number 80NSSC19K1152.

\bibliographystyle{unsrtnat}
\bibliography{Refs}

\end{document}